\documentclass[12pt]{amsart}
\usepackage[cmtip,matrix,arrow,curve]{xy}
\usepackage{graphicx}
\def \fg{{\mathfrak g}}
\def \fh{{\mathfrak h}}
\def \fk{{\mathfrak k}}

\def \ii{{\sqrt{-1}}}

\def \ft{{\mathfrak t}}

\def \fn{{\mathfrak n}}

\def \u{{\mathfrak u}}
\def \su{{\mathfrak s}{\mathfrak u}}
\def \sl{{\mathfrak s}{\mathfrak l}}

\def \C{{\mathbb C}}

\def \N{{\mathbb N}}

\newcommand\func[1]{\operatorname{\mathrm{#1}}}
\newtheorem{theorem}{Theorem}[section]

\newtheorem{proposition}[theorem]{Proposition}

\renewcommand{\t}{\mathfrak{t}}
\newcommand{\g}{\mathfrak{g}}

\newcommand{\Ad}{\func{Ad}}

\numberwithin{equation}{section}

\begin{document}

\baselineskip=16pt

\title[On sums of admissible coadjoint orbits]{On sums of admissible coadjoint orbits}

\author[A. Eshmatov]{Alimjon Eshmatov}
\author[P. Foth]{Philip Foth}

\address{Department of Mathematics, University of Arizona, Tucson, AZ 85721-0089}

\email{foth@math.arizona.edu}

\email{alimjon@math.arizona.edu}



\date{August 12, 2008}

\begin{abstract}
Given a quasi-Hermitian semisimple Lie algebra, we describe possible spectra of the sum of two
admissible elements from its dual vector space.
\end{abstract}

\maketitle
\section{Introduction}

For a compact Lie algebra $\fk$, the question of finding possible spectra of the sum of
two elements from $\fk^*$ has been answered in great details and turned out to be related to
many different areas of mathematics, including representation theory, combinatorics,
symplectic geometry, geometric invariant theory and others, see \cite{Ful} and references therein.
However, the non-compact case seems to
remain untouched, for a reason that the non-compact coadjoint orbits do not possess
the necessary nice properties and the sum of two non-compact orbits can in general cover pretty much
arbitrary spectra.

However, there is a class of coadjoint orbits in a quasi-Hermitian semisimple Lie algebra $\fg^*$,
which we call admissible, and
which share certain properties of the compact case. For example, the moment map for the maximal torus
action on such orbits is proper and the image is semi-bounded. The set of admissible orbits forms
a double cone (i.e. the union of a cone and its negative), and if we denote by $\fg^*_{\rm adm}$
the interior of one of its halves. In this setup, it now
makes sense to pose the question about the possible spectra of the sum of the orbits
of two chosen elements from the dual space of a maximal torus $\ft^*$. For example, if $a$ and $b$ are
positive real numbers and $A$ and $B$ are two matrices, which are ${\rm SU}(1,1)$-conjugate to
${\rm diag}(a, -a)$ and ${\rm diag}(b,-b)$ respectively, then possible eigenvalues $(c, -c)$ of their sum $A+B$
must necessarily satisfy the reversed triangle inequality: $c\ge a+b$.

In general, given two admissible $A, B\in\fg^*_{\rm adm}$ with prescribed spectra $\Lambda_A$ and
$\Lambda_B$ respectively, we show that the possible spectra of $A+B$ belong to a convex polyhedral set
$\left( \Pi + {\mathcal C}\right) \cap\ft^*_+$, where $\Pi$ is the polytope that solves the
problem for the maximal compact subalgebra $\fk\subset\fg$, and ${\mathcal C}$ is the cone defined by
the positive non-compact roots. The main ingredients here are the Weinstein's generalization of
the Kirwan convexity theorem for semisimple Lie groups \cite{Wein},
the Hilgert-Neeb-Plank abelian convexity theorem for non-compact manifolds \cite{HNP},
the Bates-Lerman local normal form \cite{BL}, and the Sjamaar's construction of local cones \cite{Sja}.

In the last section we describe a relationship with the tensor products of the
holomorphic discrete series representations.

\section{Admissible coadjoint orbits}

Let $G^\C$ be a complex semisimple Lie group and $\fg^\C$ its Lie algebra.
Recall the construction of real forms of $\fg^\C$ using Vogan diagrams \cite[Theorem 6.88]{Knapp}.
First, we need to fix some data for $\fg^\C$.
Choose a Cartan subalgebra $\fh$ of $\fg^\C$ and let $\Delta$ be the
root system for $(\fg^\C, \fh)$.
Fix a choice of positive roots $\Delta^+$ and let $\Sigma$ be the basis of simple roots.
Let $\ll, \, \gg$ be the Killing form of $\fg^\C$ and let root vectors $\{E_\alpha:
\alpha \in \Delta\}$ be
chosen such that $[E_\alpha, E_{-\alpha}] = H_\alpha$ for each $\alpha \in \Delta^+$,
where $H_\alpha$ is the unique element of $\fh$ defined by $\ll H,
H_\alpha\gg=\alpha(H)$ for all $H\in\fh$, and such that the numbers $m_{\alpha, \beta}$
given by $[E_\alpha, E_\beta] =  m_{\alpha, \beta}E_{\alpha+\beta}$ when $\alpha+\beta\in \Delta$ are real.
Define a compact real form $\u$ of $\fg^\C$
as
\[
\u = {\rm span}_{{\mathbb R}} \{\ii H_\alpha, \, X_\alpha:=E_\alpha-E_{-\alpha},
\, Y_\alpha:=\ii(E_{\alpha}+E_{-\alpha}) \}\ ,
\]
and let $\theta$ be the complex conjugation of $\fg$
defining $\u$.

Given a Vogan diagram $v$ for $\fg^\C$, normalized
(i.e. at most one painted root in each connected component) and with the trivial automorphism,
 let $t_v$ be the unique element in the adjoint group of $\fg^\C$ such that
$$
{\rm Ad}_{t_v}(E_\alpha)=\begin{cases}
E_\alpha \mbox{  if } \alpha \mbox{ is a blank vertex in $v$} \\
- E_\alpha \mbox{ if } \alpha \mbox{ is the painted vertex in $v$}
\end{cases}
$$
Define a complex conjugate linear involution
$$
\tau:={\rm Ad}_{t_v}\circ\theta.
$$

We use $\fg = (\fg^\C)^{\tau}$ to denote the
real form of $\fg^\C$ defined by $\tau$.
Then $\theta$ restricts to a Cartan involution of $\fg$,
and $\fh^{\tau}=\ft$
is a  compact Cartan subalgebra of $\fg$. The complexification of $\tau$ is
\begin{equation}
\label{eq_gamma_v}
\gamma := \tau \theta = \theta \tau = {\rm Ad}_{t_v} .
\end{equation}
 Since $\gamma(\Delta^+) = \Delta^+$,
the Vogan diagram of $\fg^\C$ associated to the triple $(\fg, \ft, \Delta^+)$
is $v$. Moreover, every semisimple real Lie algebra of inner type (i.e. those
that have a compact Cartan subalgebra)
can be obtained this way \cite{Knapp}.

We naturally call the roots from $\Delta^+$, as well as their negatives,
 {\it compact} if they do not have a painted root
in their decomposition into the sum of simple roots from $\Sigma$, and {\it non-compact} otherwise.
In what follows, we assume that the Lie algebra $\fg$ is quasi-Hermitian (i.e. maximal compact
subalgebras in every simple factor have non-trivial centers) and that the system of positive
non-compact roots $\Delta^+_{\rm nc}$ is adapted \cite{Neeb}, i.e. is invariant under the baby Weyl group
$W_\fk$, the Weyl group of the pair $(\fk, \ft)$, where $\fk$ is the maximal compact
subalgebra $\fk=(\fg)^\theta$.

Consider the dual vector space $\fg^*$ of the Lie algebra $\fg$, which can be identified using
the positive definite inner product $-\ll\cdot,\theta\cdot\gg$. Now we recall the definition of
two invariant cones in $\ft$:
$$
C_{\rm min} = {\rm Cone}\{ \ii [E_{\alpha}, \theta E_{\alpha}], \alpha\in\Delta^+_{\rm nc} \}
\ \ {\rm and} \ \
C_{\rm max} = \ii(\Delta^+_{\rm nc})^*\ .
$$
Now we define the open convex cone of admissible elements $\ft^*_{\rm adm}\subset\ft^*$ as
the relative interior of the dual to the minimal cone, $C^*_{\rm min}$. Using the above pairing, we can think
of $\ft^*$ as a subspace of $\fg^*$ and define the open cone of admissible elements in the latter
as $\fg^*_{\rm adm}={\rm Ad}^*(\ft^*_{\rm adm})$. (What we call here admissible is usually called strictly
admissible in the literature, but this should not lead to confusion.)

Consider $X\in\ft^*_{\rm adm}$ and its coadjoint orbit ${\mathcal O}_X$.
The co-adjoint action of $G$ on ${\mathcal O}_X$ is proper, as well as the corresponding $T$-moment map
${\mathcal O}_X\to \ft^*$, given by the projection $\fg^*\to\ft^*$ dual to the inclusion. The
image of this moment map is the sum of the polytope ${\rm Conv}(W_{\fk}.X)$ and the cone spanned by
$\ii\Delta^+_{\rm nc}$, see \cite{HNP}.
\medskip

\noindent{\it Example.} Consider $\fg^\C=\sl(n, \C)$ identified with the space of traceless
complex matrices, $\fg=\su(p,q)$ - the subspace of matrices $B$, satisfying $BJ_{pq} + J_{pq}B^*=0$ and let
$\fg^*$ be its dual vector space, which is identified with the space $\sqrt{-1}\cdot\fg$ of
pseudo-Hermitian matrices $A$, satisfying $AJ_{pq}=J_{pq}A^*$. Here
$J_{pq}={\rm diag}(\underbrace{1, ..., 1}_p, \underbrace{-1, ..., -1}_q)$ and $A^*$ is the conjugate transpose.
In the block form,
$$
A=\left(
\begin{array}{cc}
H_p & B \\
{} & {} \\
-{\bar B^T} & H_q
\end{array}
 \right)\ \ ,
$$
where $H_p$ and $H_q$ are $p\times p$ and $q\times q$ Hermitian symmetric matrices
respectively and $B$ is a complex $p\times q$ matrix. If we, as usual, take the upper-triangular matrices
as the Borel subalgebra of $\sl(n, \C)$ defined by the positive roots, then the cone of admissible elements
$\fg^*_{\rm adm}$ would consist of matrices, which are ${\rm SU}(p,q)$-conjugate to the diagonal
(and thus real) matrices of the form
${\rm diag}(\lambda_1,...,\lambda_p, \mu_1, ..., \mu_q)$ such that $\lambda_i > \mu_j$ for all pairs $i,j$. We can
certainly assume that $\lambda$'s are arranged in the non-increasing order
$\lambda_1\le\lambda_2\le \cdots \le \lambda_p$ and $\mu$'s are in the non-decreasing order
$\mu_1\ge \mu_2\ge\cdots\ge\mu_q$ (this is done for convenience), and thus the condition of admissibility becomes rather simple: $\lambda_1 > \mu_1$.

\section{Non-abelian convexity}

Let $\Lambda_A, \Lambda_B\in \ft^*_{\rm adm}$ and let
${\mathcal O}_A$ and ${\mathcal O}_B$ be the corresponding coadjoint orbits.
We would like to describe the intersection $({\mathcal O}_A + {\mathcal O}_B)\cap\ft^*_+$. In other words,
given admissible $A$ and $B$ with fixed spectra, we would like to know possible values of the spectrum
of their sum $A+B$.

Both orbits ${\mathcal O}_A$ and ${\mathcal O}_B$ are Hamiltonian $G$-spaces with the moment maps simply
given by their inclusions into $\fg^*$. Their product $({\mathcal O}_A \times {\mathcal O}_B)$
is a Hamiltonian $G$-space as well, with the moment map $\Phi$ equal to the sum of the two inclusions.
We note that, according to a generalizatoin due to Weinstein \cite{Wein}
of the Kirwan convexity theorem, the intersection $\Phi({\mathcal O}_A \times {\mathcal O}_B)\cap\ft^*_+$
is a convex polyhedral set ${\mathcal S}_{AB}$, describing exactly all possible spectra of such sums $A+B$.
In this section we will give a more detailed description of this set.

We will start by describing the local convexity data. For brevity, denote
$M={\mathcal O}_A \times {\mathcal O}_B$ and $\omega$ the product symplectic form on $M$.
For any subset $Q\subset M$ denote by
${\mathcal S}(Q)$ the image of $\Phi(Q)$ in $\ft^*_{\rm adm}=\fg^*_{\rm adm}/G\subset \ft^*_+$.
For a given point $m\in M$, let $G_m$ be its stabilizer, and let $G_y$ be the stabilizer
of $y=\Phi(m)$ with respect to the
co-adjoint action of $G$. Note that since $\Phi(M)$ is in the admissible cone, where the action of $G$ is proper,
$G_y$, and its subgroup $G_m$, are both compact.
Let $\fg_m$ and $\fg_y$ denote the corresponding subalgebras of $\fg$. Using the Killing form as before,
we can also think of $\fg^*_m$ as a subspace of $\fg^*_y$, and the latter as a subspace of $\fg^*$.
Let $\fg_m^\perp$ be the annihilator of $\fg_m$ in $\fg^*_y$ and $\fg_y^\perp$ be the annihilator of $\fg_y$ in $\fg$.
Then we have an $G_m$-equivariant splitting $\fg^* = \fg^* _y \oplus \fg_m^\perp \oplus \fg_y^\perp$and we denote
$i \, : \, \fg^*_m \rightarrow \fg^* _y$ and $j \, : \, \fg^*_y \rightarrow \fg^*$ the corresponding injections.
Let also ${\mathcal O}_m$ be the $G$-orbit through $m$ and ${\mathcal O}_y$ - the coadjoint orbit of $y$.
Consider the symplectic vector space
\begin{equation}
V=(T_m{\mathcal O}_m)^\perp/((T_m{\mathcal O}_m)^\perp\cap T_m{\mathcal O}_m),
\label{eq:v}
\end{equation}
where $\perp$ stands for the symplectic perp. This space has a natural linear symplectic action of $G_m$,
with moment map $\Psi_{V}: V\to\fg^*_m$. A theorem of Bates and Lerman \cite{BL}, extending the results
of Guillenin-Sternberg and Marle, asserts that:

\begin{proposition}There exists a $G$-invariant neighbourhood $U$ of ${\mathcal O}_m$ in $M$ and a
$G$-invariant neighbourhood $U_0$  of the zero section of the vector bundle
$G\times_{G_m}(\fg_m^\perp\times V)\to G/G_m$, and a $G$-equivariant symplectomorphism
$\eta: U_0\to U$ such that
$$
\Phi(\eta(g, X, v))={\rm Ad}^*_g(y+j(X+i(\Psi(v))))\ .
$$
\end{proposition}

Next, recall the constructive proof of the non-abelian Convexity Theorem, due to Sjamaar, as explained
in \cite{GSj}. This result gives a concrete description of the local structure of the Kirwan polytope,
and the convexity is explained in terms of ``$T$ to $K$ induction''.  Denote by $\Pi$
a convex polytope in $\ft^*_+$, which appears in the Kirwan convexity theorem for the
maximal compact subgroup $K\subset G$:
$$
\Pi := {\mathcal S}\left( \Phi({\rm Ad}^*_K(\Lambda_A)\times {\rm Ad}^*_K(\Lambda_B))  \right) .
$$

Sjamaar's proof of the Kirwan convexity theorem readily extends to
our current setup, mainly due to the fact that the actions are proper and all the
stabilizers are compact.

\begin{theorem} The image ${\mathcal S}_{AB}={\mathcal S}(M)$ in $\ft^*_{\rm adm}$ is the
intersection of local moment cones and is given by a convex polyhedral set
$(\Pi + {\rm Cone}(\ii\Delta^+_{\rm nc}))\cap\ft^*_+$. A point $y$ is
an extremal points of ${\mathcal S}(M)$ if and only if $\fg_y=[\fg_y,\fg_y] +\fg_m$,
where $m\in\Phi^{-1}(y)$.
\end{theorem}

\noindent
\begin{proof}
In order to prove this fact we need two ingredients: local convexity theorem and connectivity of fibers.
The latter condition is guaranteed by Weinstein's result \cite[Theorem 3.3]{Wein}. On the other hand, by Proposition 3.1
in a $G$-invariant neighborhood of point $m$ we have a simple canonical model given by $Y\, = \, G \times ^{G_m}
(\mathfrak{g}^{\bot}_m \times V)$. The symplectic manifold $Y$ can be realized as a symplectic quotient of $X=G \times
 \g^* _{y} \times T^* \, G_{y} \times V$ by the $G_{y} \times G_{m}$ action. Now, combining this with
 the commutativity of reduction in stages, we are able to compute the local moment cones.

Let us first describe the Hamiltonian structure on $X$ more explicitly. The manifold $G\times \g^*_{y}$ carries a
natural closed two-form (the minimal-coupling form, see \cite{GS}), which is non-degenerate in a $G$-invariant neighborhood
of $G\times \{0\}$. Now, identifying $T^* \, G_{y}$ with $G_{y} \times \g^* _y$ by means
of left translations, the action of $G_y \times G_m$ on $X$ is given by:
$$
(g,h) \, . \, (l , \xi , k, \mu , v) = (lg^{-1} , \Ad^*_g \xi , g k h^{-1}, \Ad^*_h \mu , hv)\ ,
$$
which is Hamiltonian with the moment map:
$$
\Psi (l , \xi , k, \mu , v) =(-\xi + \Ad^*_{k} \mu, -\pi(\mu) + \Phi_V (v))\ ,
$$
where $\Phi_V$ is the moment map for the linear symplectic action of $G_m$ on $V$ and
$\pi: \g^*_y \rightarrow \g^*_m$ is the natural projection.
Note, since both $G_y$ and $G_m$
are compact, $\Psi$ is proper. Since $(-y,0)$ is a regular value and $G_y \times G_m$ acts freely,
the symplectic quotient $X_{y} = \Psi^{-1} (-y,0)/(G_y \times G_m)$
is a smooth manifold.
The space $X_y$ carries a natural $G$-action inherited from $X$, which is Hamiltonian. Now, in order to
obtain the desired description of $Y$, consider the map $\phi$ from
$G \times \mathfrak{g}^{\bot}_m \times V$ to $\Psi^{-1} (-y,0)$ defined by:
$$
\phi(l, \mu, v) \, = \, (l, \mu+ \Phi_V (v)+ y, 1, \mu+ \Phi_V (v), v)\ ,
$$
which clearly is $G$-equivariant. Therefore $\phi$ descends to the following $G$-equivariant diffeomorphism:
$$
\bar{\phi}: Y \rightarrow X_{y}\ .
$$
Now, using the equivariant Darboux theorem, one can see that this map is in fact a symplectomorphism.
If we perform the above reduction in stages, namely first
with respect to $G_m$, and then by $G_y$, we can present $Y$ as an iterated bundle:
$$
Y \cong G \times^{G_y} (G_y \times^{G_m} (\mathfrak{g}^{\bot}_m \times V)) = G \times^{G_y} \tilde Y
$$
over the coadjoint orbit, diffeomorphic to $G/G_{y}$, with the fiber
$\tilde Y=G_y \times^{G_m} (\mathfrak{g}^{\bot}_m \times V)$. The space $\tilde Y$ is a Hamiltonian $G_y$-space with
the moment map $\Phi|_{\tilde Y}$, which is a restriction of $\Phi$ to $\tilde Y$.
Therefore we can write $\Phi$ as the composition of maps:
$$
\xymatrix{
Y \, = \, G \times^{G_y} \tilde Y  \ar[r]^{\ \ id \times \Phi} &  G \times^{G_y} \g^* _{y} \ar[r]^{\ \ \iota} & \g^*
}
$$
where $\iota([g, \xi])= \Ad^*_g (\xi)$. On the other hand,
$\g^*_{y}$ is a slice at $y$ for a coadjoint action on $\g^*$.
Thus, restricting $\iota$ to a sufficiently small $G$-invariant neighborhood of $[g, \xi]$,
we obtain a $G$-equivariant
embedding into an open neighborhood of $y$. Hence there is a $G$-invariant
neighborhood $U$ of $[1,0,0]$ in $Y$, so that
$\Phi$ becomes the bundle map of associated bundles over $G/G_y$.
Therefore $U \cap \Phi^{-1}(\g^*_y)=U \cap Y$, and the
image $\Phi(U)$ is a bundle over $G.y$ with the fiber $\Phi(U \cap \tilde Y)$.
If $U$ is small enough, we get:
$$
\mathcal{S}(U) = \Phi(U) \cap \mathfrak{t^*_{+}}= \Phi(U \cap Y) \cap \mathfrak{t}^* _{+}=
\mathcal{S}(U \cap \tilde Y)\ ,
$$
where $ \mathcal{S}(U \cap \tilde Y)$ is the
moment image of $U \cap \tilde Y$ as a Hamiltonian $G_y$-space. Hence the description
of local moment map reduces to computing the local moment map of the Hamiltonian $G_y$-space $\tilde Y$.
Note that since $G_y$ compact
and $\Phi(U \cap \tilde Y)$ is proper,
we can use the vertex criterion given in \cite[Theorem 6.7]{Sja}, which implies that if $y$ is
a vertex, then $\g_y = [\g_y , \g_y] + \g_m$, or, equivalently, $G_y = [G_y , G_y] G_m$.
In particular, if $y$ is an interior point
of $\t^* _+$ and $T$ fixes $m$, then
$m \in (W_{\fk}.A , W_{\fk}. B)$.
%  I have reservations about that last sentence.

Now, in what follows, we give an explicit description of the local moment cone.
The space $\tilde Y$ is a symplectic quotient of
$T^* G_y \times V \cong G_y \times \g^*_y \times V$ by $G_m$ with the moment map:
$$
\widetilde \Psi (g, \xi, v) = -\pi(\xi) + \Phi_V(v) + \pi(y)\ ,
$$
where we shifted the moment map by $y$. Let us assume that $G_m$ is abelian, which is the case 
for an open dense
set of elements in $M$. This follows from the fact that the isotropy group $G_y$, which contains $G_m$,
is a subgroup of $T$ for a dense open subset of elements  $y \in \ft^*_{+}$   .
Therefore we can think of $G_m$ as a subgroup of $T$. The $G_m$-moment map image of $V$ is the cone:
$$
\mathcal{C}_m=\{ \sum^n_{i=1} t_i \alpha_i \, |  \, t_i \geq 0 \}\ ,
$$
 where $\{\alpha_i\}_{i=\overline{1,n}}$ are the weights of the
 representation of $G_m$ on $V$. For $\lambda \in \mathfrak{t^*_{+} }$, consider the coadjoint
 orbit $G_y . (\lambda) $ through $\lambda$.
 If we regard $T^* G_y \times V$ as a Hamiltonian $G_y$-space and reduce it with
 respect to $G_y$, we obtain $ -\, G_y . (\lambda) \times V$. It is a Hamiltonian
$G_m$-space with the moment map image:
$$
- \pi(\mathtt{Conv} (W_y . \lambda))+\mathcal{C}_m + \pi(y).
$$
Thus the reduction of $ -\, G_y . (\lambda) \times V$ with respect to $0 \in \g^* _m$ is non-empty 
if and only if:
$$
0 \in - \pi(\mathtt{Conv} (W_y . \lambda))+\mathcal{C}_m + \pi(y).
$$
If we regard $T^* G_y \times V$ as a $G_m$-space and first reduce with respect to $0$ to obtain
$\widetilde Y$, and then reduce with respect to $G_y.(\lambda)$, then
we obtain the very same space, according to the general result
for reduction in stages \cite{OR}. Hence this space is non-empty
if and only if $G_y. (\lambda)$ is in the moment image of $\widetilde Y$.
This implies that the local moment cone is the set:
$$
\{ \, \lambda \in \mathfrak{t}^* _{+} \, | \, \pi(\mathtt{Conv} (W_y . \lambda)) \,
\cap (\mathcal{C}_m + \pi(y)) \neq \emptyset\}.
$$
or, equivalently, that there is some neighborhood $U$ of $y$ in $\ft^*_+$ such that:
$$
U \cap\mathcal{S}(\tilde Y) = U \cap (\pi^{-1} _{\ft} (\mathcal{C}_m )+y )\ ,
$$
where $\pi_{\ft} : \t^* \rightarrow \fg_m ^*$ is the natural projection.
Combining this with the previously mentioned theorem of Weinstein \cite{Wein},
which also asserts the connectedness of the fibers of ${\mathcal S}: M\to\ft^*_{\rm adm}$,
we conclude that ${\mathcal S}(M)$ is closed convex polyhedral subset of $\ft^* _{\rm adm}$.

In order to obtain the description  stated in the theorem,
we recall a simple result about closed convex sets \cite[Proposition 5]{HNP}.
Let $\mathcal{C}$ be a closed convex set in a vector space $V$. Then:
$$
\mathcal{C}= \mathtt{Conv}({\rm Ext}(\mathcal{C})) + \lim (\mathcal{C})\ ,
$$
where ${\rm Ext}(\mathcal{C})$ is the set of extremal points of the cone, and
$\lim(\mathcal{C}): = \{ v \in V : \mathcal{C}+ v \subseteq \mathcal{C}\}$.
Now using above vertex criterion ${\rm Ext}(\mathcal{S(M)}) \, \cap \,
(\ft^*_{+})^{\circ} \, =  \, \Phi(W_{\fk}.A , W_{\fk}. B) \cap (\ft^*_{+})^{\circ}$
where $(\ft^*_{+})^{\circ}$ is the interior of $\ft^*_{+}$, we deduce that
$\mathtt{Conv}({\rm Ext}(\mathcal{S(M)}))\, \cap \,
(\ft^*_{+})^{\circ}  \, = \, \Pi \cap (\ft^*_{+})^{\circ} $. Hence we have:
$$
S(M)  \cap \, (\ft^*_{+})^{\circ}  \, = \, ( \Pi + \lim S(M) ) \cap \, (\ft^*_{+})^{\circ}\ .
$$
On the other hand, the cones at the extremal points have a very simple form.
Using \cite[Remark 5.18]{HNP} for the moment map image of coadjoint orbits,
the cone at those points is given by ${\rm Cone}(\ii\Delta^+_{\rm nc})$. Therefore we have:
$$
S(M)  \cap \, (\ft^*_{+})^{\circ}  \, = \,  (\Pi + {\rm Cone}(\ii\Delta^+_{\rm nc})) \cap (\ft^*_{+})^{\circ}\ .
$$
And by continuity we conclude:
$$
S(M)\, = \,  (\Pi + {\rm Cone}(\ii\Delta^+_{\rm nc})) \cap \ft^*_{+}\ .
$$
\end{proof}

\medskip

\noindent{\it Example 1.} Let $G={\rm SU(2,1)}$. Identify, as before, $\fg^*$ with the space of pseudo-Hermitian
matrices of signature $(2,1)$ and take $\Lambda_A={\rm diag}(4, 1, -5)$ and $\Lambda_B={\rm diag}(2, 1, -3)$.
Then the possible eigenvalues $(\lambda_1, \lambda_2, \mu)$ of $A+B$, taken in the order
$\lambda_1 \ge \lambda_2 > \mu$, are given by:
$$
\lambda_1\ge 5, \ \ \lambda_2\ge 2, \ \ \lambda_1+\lambda_2\ge 8, \ \ {\rm and} \ \ \mu=-\lambda_1-\lambda_2.
$$
\begin{figure}
    \begin{center}
    \includegraphics[height=3in,width=3in,angle=0]{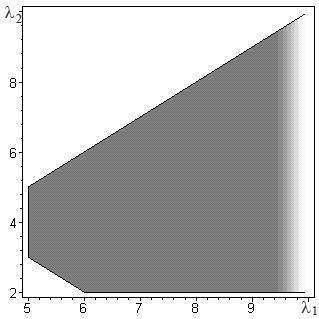}
  \end{center}
  \caption{Moment polyhedron for Example 1.}
\end{figure}

\medskip

\newpage
\noindent{\it Example 2.} Let $G={\rm SU(2,2)}$. 
Take $\Lambda_A={\rm diag}(4, 2, 1, -7)$ and $\Lambda_B={\rm diag}(3, 2, 1, -6)$.
Then the possible eigenvalues $(\lambda_1, \lambda_2, \mu_1 ,\mu_2)$ of $A+B$, taken in the order:
$$
\lambda_1 \ge \lambda_2 > \mu_1 \geq \mu_2 \ ,
$$
are given by:
$$
\lambda_1\ge 6, \ \ \lambda_2\ge 4, \ \ \lambda_1+\lambda_2\ge 11, \ \  \lambda_1+\lambda_2 + \mu_1 \ge 6, \ \ \mu_1 \leq 2 \ \ {\rm and} \ \ \mu_2=-\lambda_1-\lambda_2-\mu_1 .
$$
If we project this region onto the first three coordinates, then we will obtain a polyhedron sketched in Figure 2. 

\begin{figure}
    \begin{center}
    \includegraphics[height=5in,width=5in,angle=0]{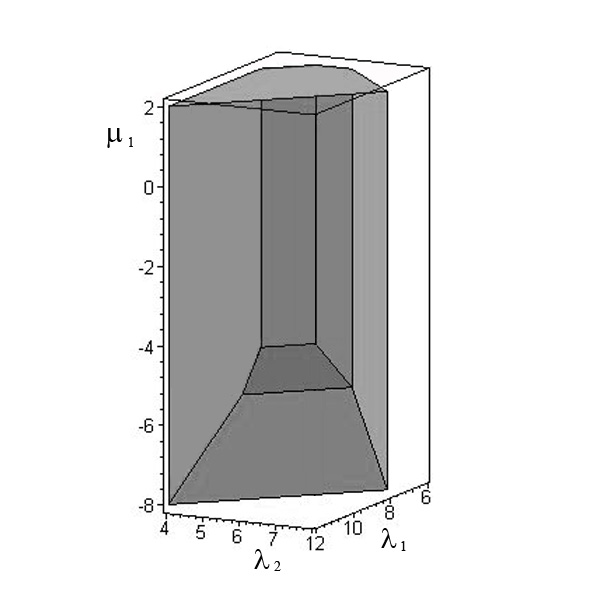}
  \end{center}
  \caption{Moment polyhedron for Example 2.}
\end{figure}
\section{Relationship with representation theory}

Let $\Lambda\in\ft^*_{\rm adm}$ be a dominant integral weight with respect to the compact positive roots and
let $V$ be a irreducible unitary $K$-module with highest weight $\Lambda$.
Following Harish-Chandra, one can construct a unique unitary irreducible representation $\rho_\Lambda$
of $G$, such that the corresponding representation of $G^\C$ has highest weight $\Lambda$. The underlying space of $\rho_\Lambda$ is $V\otimes {\mathcal U}_+$, where ${\mathcal U}_+$ is the
universal enveloping algebra of the nilradical $\fn_+$ spanned over $\C$ by the non-compact positive roots.
(Note that we have a bit different convention and take $\fn_+$ instead of $\fn_-$, but also
our $\Lambda$ is in $\ft^*_{\rm adm}$ and not in its negative.) Such a $\rho_\Lambda$ is called a holomorphic discrete series reresentation and is a generalized Verma module as a representation of ${\mathcal U}_{\fg}$.

The weights of $\rho_\Lambda$ have the form $\lambda+\N .\Delta^+_{\rm nc}$, where $\lambda$ is
a weight of $V$ and $\N$ is the semigroup of non-negative integers.

Moreover, a theorem of Repka \cite{Repka}
says that for two such representations $\rho_A$ and $\rho_B$, their tensor product decomposes into
the sum of subspaces, whose weights are sums of weights of $\rho_A$ and $\rho_B$, and all are of the form
$ \lambda_{AB} + \N .\Delta^+_{\rm nc}$, where $\lambda_{AB}$ is a weight of $V_A\otimes V_B$.

Thus one concludes that the tensor product of holomorphic discrete series representations
of $G$ with highest weights $\Lambda_A$ and $\Lambda_B$ decompose into the direct sum of representations
with highest weights given by the lattice points in the convex polyhedral set
${\mathcal S}({\mathcal O}_A\times{\mathcal O}_B)$, with finite multiplicities.
It would be interesting to find a direct approach to proving this result, as well as to establish that
the (finite) multiplicities of the representations appearing as summands in such tensor products
correspond to counting the lattice points in certain polyhedral sets, similar to the compact case.

\section*{Acknowledgements}

We would like to thank Reyer Sjamaar for reading a preliminary version of this manuscript
and making a number of useful comments. We are
also grateful to Alexander Dvorsky and Sam Evens for useful communications.

%%%%%%%%%%%%%%%%%%%%%%%%%%%%%%%%%%%%%%%%%%%%%%%%%%%%%%%%%%%%%%%%

\end{document}